# THERE ARE NO HERONIAN TRIANGLES WITH THREE INTEGER MEDIANS

Logman Shihaliev logman1@list.ru +994505149553

ORCID: 0000-0003-1063-4712

**ABSTRACT**. The relevance of this paper lies in the fact that it resolves two previously unsolved open problems. In the first part of the paper, a new lemma is proved, from which it follows that if there exists a triangle with integer sides and medians, then there necessarily exists another triangle, not similar to it, with the same properties. In other words, such triangles can exist only in pairs.

In the second part of the paper, by transforming known formulas, a new theorem is established in the form of a universal identity valid for all triangles. The focus of this theorem is the proof of the nonexistence of Heronian triangles with three integer medians.

We arrive at the conclusion that, among the seven elements of a triangle (three sides, three medians, and the area), only six can be integers.

It should be noted that if the above universal identity is considered as a Diophantine equation, then such a Diophantine equation does have solutions. As examples, one may consider the following pair of triangles. A triangle with sides $146, 102, 52$, area $1680$, and two integer medians $35$ and $97$. Also, according to the conditions of the above lemma, there exists another triangle with medians $219, 153, 78$, area $5040$, and two integer sides $70$ and $194$.

This version of the article is different from the previous ones because everything has been simplified to the level of a straight-A high school student.

Another noteworthy result is presented in **Corollary 2**, where the paradox of "twin" triangles is identified.

## PART I

**LEMMA.** If there exists a triangle with rational sides and rational medians, then there necessarily exists another triangle, not similar to it, with rational sides and rational medians.

**PROOF OF THE LEMMA.**

Consider a triangle $\Delta ABC$ with rational sides and rational medians (Fig. 1). Using the triangle $\Delta ABC$, we construct a triangle $\Delta OPC$.

To this end, starting from point $C$, we draw the segment $CP \| OB$ until it intersects the extension of the median $AF = m_a$ at point $P$.

As a result, the triangles $\Delta APC$ and $\Delta AOE$ are similar. We write the ratio:

$CP{:}\,EO = AP{:}\,AO = AC{:}\,AE = 2{:}\,1$

Taking into account the properties of the medians of the triangle $(\Delta ABC)$ with respect to the sides of the triangle $\Delta OPC$, we obtain

$$CP = OB = \frac{2}{3}m_b, \quad OC = \frac{2}{3}m_c, \quad OP = OA = \frac{2}{3}m_a \qquad \textbf{(1)}$$

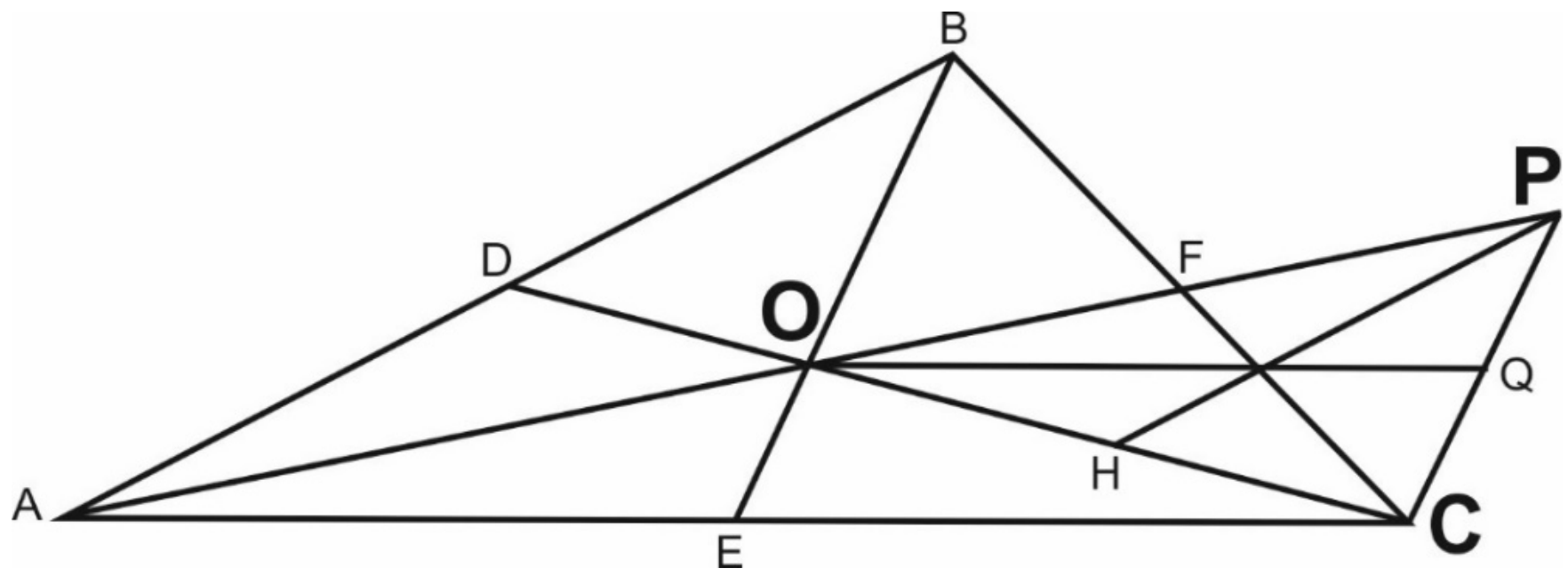

**Fig1**

Here (**Fig. 1**): $AB = c, BC = a, \; AC = b, AF = m_a, BE = m_b, \; CD = m_c,$
$OP = \frac{2}{3}m_a, CP \| OB, CP = OB = \frac{2}{3}m_b, OC = \frac{2}{3}m_c \;\; OQ = EC = \frac{1}{2}b,$
$HP = DB = \frac{1}{2}c, CF = \frac{1}{2}a, A_{\Delta ABC} \; = 3A_{\Delta OPC}$

And for the medians of the $\Delta OPC$, we get:

$$FC = \frac{1}{2}BC = \frac{1}{2}a, \quad OQ = \frac{1}{2}AC = \frac{1}{2}b, \quad HP = \frac{1}{2}AB = \frac{1}{2}c \qquad \textbf{(2)}$$

In other words, the triangle $\Delta OPC$ is constructed by parallel translations of the $\frac{2}{3}$ -parts of the medians of the triangle \triangle ABC. The medians of the triangle $\Delta OPC$, in turn, are constructed by parallel translations of the $\frac{1}{2}$ -parts of the sides of the triangle $\Delta ABC$. This implies that all sides and medians of the triangle $\Delta OPC$ are also rational.

Since $A_{\Delta OFC} = \frac{A_{\Delta ABC}}{6}$, we get $A_{\Delta OPC} \; = 2 \cdot A_{\Delta OFC} = \frac{A_{\Delta ABC}}{3}$. We write:

$$\frac{A_{\Delta OPC}}{A_{\Delta ABC}} = \frac{1}{3} \qquad \textbf{(3)}$$

The ratio of the areas of similar triangles with rational sides is equal to the square of a rational similarity coefficient. But the fraction $\frac{1}{3}$ is not the square of a rational number.

**THE LEMMA IS PROVED.**

**COROLLARY 1**: If (1) and (2) are applied to triangle $\Delta OPC$, we obtain a triangle similar to $\Delta ABC$.

For example, let the sides of triangle $\Delta ABC$ be $136, 170, 174$, and let its medians be $158, 131, 127$. Then we obtain triangle $\Delta OPC$ with integer sides $316, 262, 254$, and integer medians $204, 255, 261$. Continuing this process, we obtain the next triangle with sides $204, 255, 261$, and medians $237,\ \ 196.5,\ \ 190.5$.:

$$\frac{204}{136} = \frac{255}{170} = \frac{261}{174}.$$

Let's write down the well-known mathematical identities (as a reference):

$$\begin{cases} a = \frac{2}{3}\sqrt{2m_b^2 + 2m_c^2 - m_a^2} \\ b = \frac{2}{3}\sqrt{2m_a^2 + 2m_c^2 - m_b^2} \\ c = \frac{2}{3}\sqrt{2m_a^2 + 2m_b^2 - m_c^2} \end{cases} \qquad \begin{cases} m_a = \frac{1}{2}\sqrt{2b^2 + 2c^2 - a^2} \\ m_b = \frac{1}{2}\sqrt{2a^2 + 2c^2 - b^2} \\ m_c = \frac{1}{2}\sqrt{2a^2 + 2b^2 - c^2} \end{cases} \qquad \textbf{(4)}$$

## PART II

**THEOREM**. For any triangle, and for any ordering of its sides (there are 3! = 6 possible permutations), the following identity holds:

$$\left(\left(\frac{a}{2}\right)^2 - \left(\frac{b}{2}\right)^2\right)^2 + (A_{\Delta ABC})^2 = \left(\frac{c}{2}m_c\right)^2. \qquad \textbf{(5)}$$

**PROOF OF THE THEOREM.**

For the area of the triangle ΔABC, we write down Heron's formula

$$A_{\Delta ABC} = \frac{1}{4}\sqrt{(a+b+c)(a+b-c)(a-b+c)(-a+b+c)} \qquad \textbf{(6)}$$

Note that Heron's formula (6) is provable using trigonometry (without applying the Pythagorean theorem).

Let us transform (6). The purpose of this transformation is to explicitly incorporate another element of the triangle (in this case, a median) into formula (6).

We know that a triangle with given side lengths is unique; therefore, its area is determined, and these four elements $(a, b, c, A_{\Delta ABC})$ are present in (6). The presence of a fifth element of the triangle (the fifth is not superfluous) in a single formula creates the possibility for more efficient analysis.

We rewrite (6) as follows:

$$A_{\Delta ABC} = \frac{1}{4}\sqrt{2(a^2b^2 + a^2c^2 + b^2c^2) - (a^4 + b^4 + c^4)}$$

$$16A_{\Delta ABC}^2 = -a^4 - b^4 - c^4 + 2a^2b^2 + 2a^2c^2 + 2b^2c^2$$

$$a^4 - 2a^2b^2 + b^4 + 16A_{\Delta ABC}^2 = 2a^2c^2 + 2b^2c^2 - c^4$$

$$\frac{a^4 - 2a^2b^2 + b^4}{16} + A_{\Delta ABC}^2 = \frac{2a^2c^2 + 2b^2c^2 - c^4}{16}$$

$$\frac{a^4-2a^2b^2+b^4}{16}+A_{\Delta ABC}^2=\frac{c^2(2a^2+2b^2-c^2)}{4\cdot 4}\Rightarrow\frac{(a^2-b^2)^2}{16}+A_{\Delta ABC}^2=\frac{c^2\cdot m_c^2}{4}$$

$$\left(\left(\frac{a}{2}\right)^2-\left(\frac{b}{2}\right)^2\right)^2+A_{\Delta ABC}^2=\left(\frac{c}{2}m_c\right)^2 \qquad \textbf{(7)}$$

Since equality (7) is symmetric with respect to the sides of the triangle, we obtain six ($3! = 6$) equivalent identities:

$$\begin{cases}\left(\pm\left(\frac{a}{2}\right)^2\mp\left(\frac{b}{2}\right)^2\right)^2+(A_{\Delta ABC})^2=\left(\frac{c}{2}m_c\right)^2\\ \left(\pm\left(\frac{a}{2}\right)^2\mp\left(\frac{c}{2}\right)^2\right)^2+(A_{\Delta ABC})^2=\left(\frac{b}{2}m_b\right)^2\\ \left(\pm\left(\frac{b}{2}\right)^2\mp\left(\frac{c}{2}\right)^2\right)^2+(A_{\Delta ABC})^2=\left(\frac{a}{2}m_a\right)^2\end{cases} \qquad \textbf{(8)}$$

**THE THEOREM IS PROVED**.

The analysis of identities (8) will allow us to solve the following open problem.

**PROBLEM.** Are there triangles with integer sides, medians, and area? [1, 2, 3, 4, 5, 6, 7]

**SOLUTION OF THE PROBLEM**.

Suppose there exists a triangle in which all medians, sides, and the area are integers (we will call such a triangle an *integer triangle*) and that

$$(a,b,c,m_a,m_b,m_c,A_{\Delta ABC})=1 \qquad \textbf{(9)}$$

Let us write down one of the identities from (8):

$$\left(\left(\frac{a}{2}\right)^2-\left(\frac{b}{2}\right)^2\right)^2+(A_{\Delta ABC})^2=\left(\frac{c}{2}m_c\right)^2 \qquad \textbf{(10)}$$

According to the conditions in (4), it is evident that all sides of *integer triangles* (in this case, triangle $\Delta ABC$) are even. We write:

$$(a,b,c)=2^k \qquad \textbf{(11)}$$

Consequently, the numbers $\frac{a}{2},\frac{b}{2}$ and $\frac{c}{2}$ are integers.

Here $k=1$. If $k>1$, then we get:

$$(m_a,m_b,m_c)=2^p. \qquad \textbf{(12)}$$

This contradicts condition (9), since the area of Heronian triangles is an even number. Instead of (9), we obtain:

$$(a,b,c,m_a,m_b,m_c,A_{\Delta ABC})=2t>1 \qquad \textbf{(13)}$$

Here $k, p$ and $t$ are natural numbers.

**AUXILIARY LEMMA**. Exactly one of the sides of an *integer trangle* is divisible by 4.

**PROOF OF THE AUXILIARY LEMMA.**

The contradiction between (9) and (13) shows that all four sides of an *integer triangle* cannot be divisible by 4 simultaneously.

**Suppose** that none of the sides of an *integer triangle* is divisible by 4. In this case, we obtain:

$$\frac{c}{2} \equiv 1(mod2) \quad \textbf{(14)}$$

$$\left(\frac{a}{2}\right)^2 - \left(\frac{b}{2}\right)^2 \equiv 0(mod2) \quad \textbf{(15)}$$

Since the area $A_{\Delta ABC}$ is an even number, it follows that the left-hand side of identity (10) is even. Consequently, the right-hand side of (10) must also be even. This is possible only if

$$m_c \equiv 0(mod2). \quad \textbf{(16)}$$

This means that, according to formula (4), side c is divisible by 4, which contradicts our assumption (14).

**Now suppose** that two sides of the *integer triangle* (let these be $a$ and $b$) are divisible by 4, and for side $c$ we have $\frac{c}{2} \equiv 1(mod2)$. In accordance with the conditions of formula (4), this means that

$$m_a \equiv 0(mod2),\ m_b \equiv 0(mod2). \quad \textbf{(17)}$$

Taking into account that

$$\sqrt{2a^2 + 2b^2 - c^2} \not\equiv 0(mod4) \quad \textbf{(18)}$$

We obtain

$$m_c \equiv 1(mod2) \quad \textbf{(19)}$$

Hence,

$$\left(\frac{c}{2}m_c\right)^2 \equiv 1(mod2). \quad \textbf{(20)}$$

Therefore, the left-hand side and the right-hand side of equality (10) have different parity, and hence our assumption is incorrect.

**THE AUXILIARY LEMMA IS PROVED.**

Is the Pythagorean triple (10) primitive? To determine this, let us consider all possible values of $q$:

$$\left(\left(\frac{a}{2}\right)^2 - \left(\frac{b}{2}\right)^2, A_{\Delta ABC}, \frac{c}{2}m_c\right) = q \quad \textbf{(21)}$$

Let us consider condition

$$(q, 6) = 1 \text{ for } q > 1 \quad \textbf{(22)}$$

If all the sides of the triangle are divisible by $q$, then, in view of identities (4), all the medians are also divisible by $q$. Consequently, the area of the triangle is also divisible by $q$:

$$(a, b, c, m_a, m_b, m_c, A_{\Delta ABC}) = q \qquad \textbf{(23)}$$

This contradicts condition (9).

If only two sides (let them be $a$ and $b$, for which we choose identity (10) from system (9)) of the triangle are divisible by $q$, then, in view of identities (4), it follows that $m_a m_b m_c \not\equiv 0(mod q)$. Consequently, for the right-hand side of (10) we get

$$\left(\frac{c}{2}m_c\right)^2 \not\equiv 0(mod q) \qquad \textbf{(24)}$$

Result (24) contradicts condition (21).

**Remark 1**. Results (21), (23), and (24) imply that if an *integer triangle* exists, then under condition (22) the Pythagorean triple $\left(\frac{a}{2}\right)^2 - \left(\frac{b}{2}\right)^2, A_{\Delta ABC}, \frac{c}{2}m_c$ is primitive.

**Remark 2**. Using the results of the auxiliary lemma, we conclude that only one of the following three expressions is even:

$$\left(\frac{a}{2}\right)^2 - \left(\frac{b}{2}\right)^2, \quad \left(\frac{a}{2}\right)^2 - \left(\frac{c}{2}\right)^2, \quad \left(\frac{b}{2}\right)^2 - \left(\frac{c}{2}\right)^2$$

Let us write:

$$\left(\frac{a}{2}\right)^2 - \left(\frac{b}{2}\right)^2 \equiv 1(mod2), \quad \left(\frac{a}{2}\right)^2 - \left(\frac{c}{2}\right)^2 \equiv 1(mod2). \qquad \textbf{(25)}$$

**Remark 3**. In identities (4), the fractions $\frac{2}{3}$ and $\frac{1}{2}$ appear. For this reason, it is necessary to analyze the values $q = \{1, 2, 3, 6\}$ and all other values of $q$ separately.

**Now**, in (21), let us consider the possible values $q = \{1, 2, 3, 6\}$.

**If** $q = 1$, then the Pythagorean triple $\left(\frac{a}{2}\right)^2 - \left(\frac{b}{2}\right)^2, A_{\Delta ABC}, \frac{c}{2}m_c$ is primitive. It is known that, according to Euclid's formula, every primitive Pythagorean triple arises from a unique pair of coprime numbers.

Consider two cases:

1. $\left(\frac{a}{2}\right)^2 - \left(\frac{b}{2}\right)^2 = \left(\frac{a}{2}\right)^2 - \left(\frac{b}{2}\right)^2$
2. $\left(\frac{a}{2}\right)^2 - \left(\frac{b}{2}\right)^2 = w^2 - r^2 = \left(\frac{2w}{2}\right)^2 - \left(\frac{2r}{2}\right)^2$

**For the first case**, we write:

$$\begin{cases} \left(\frac{a}{2}\right)^2 - \left(\frac{b}{2}\right)^2 = \left(\frac{a}{2}\right)^2 - \left(\frac{b}{2}\right)^2 \\ A_{\Delta ABC} = 2 \cdot \frac{a}{2} \cdot \frac{b}{2} = \frac{ab}{2} \\ \frac{c}{2}m_c = \left(\frac{a}{2}\right)^2 + \left(\frac{b}{2}\right)^2 \end{cases} \qquad \textbf{(26)}$$

Unlike (26), another formula for the area $A_{\Delta ABC}$ is known:

$$A_{\Delta ABC} = \frac{1}{2}ab \cdot sin\angle C \tag{27}$$

Thus,

$$\frac{1}{2}ab \cdot sin\angle C = \frac{ab}{2} \Rightarrow sin\angle C = 1 \Rightarrow \angle C = \frac{\pi}{2} \tag{28}$$

According to (25), we choose identity

$$\left(\left(\frac{a}{2}\right)^2 - \left(\frac{c}{2}\right)^2\right)^2 + (A_{\Delta ABC})^2 = \left(\frac{b}{2}m_b\right)^2 \tag{29}$$

We write

$$\begin{cases} \left(\frac{a}{2}\right)^2 - \left(\frac{c}{2}\right)^2 = \left(\frac{a}{2}\right)^2 - \left(\frac{c}{2}\right)^2 \\ A_{\Delta ABC} = 2 \cdot \frac{a}{2} \cdot \frac{c}{2} = \frac{ac}{2} \\ \frac{b}{2}m_b = \left(\frac{a}{2}\right)^2 + \left(\frac{c}{2}\right)^2 \end{cases} \tag{30}$$

We obtain $\angle B = \frac{\pi}{2}$. **(31)**

In a single triangle, two angles cannot both be equal to $\frac{\pi}{2}$.

**For the second case**, we write:

$$\begin{cases} \left(\frac{a}{2}\right)^2 - \left(\frac{b}{2}\right)^2 = \left(\frac{2w}{2}\right)^2 - \left(\frac{2r}{2}\right)^2 \\ A_{\Delta ABC} = 2 \cdot \frac{2w}{2} \cdot \frac{2r}{2} = \frac{2w \cdot 2r}{2} \\ \frac{c}{2}m_c = \left(\frac{2w}{2}\right)^2 + \left(\frac{2r}{2}\right)^2 \end{cases}$$

Here we obtain a right triangle with two sides $(2w)$ and $(2r)$ and with the same area $A_{\Delta ABC}$. Thus, we arrive again at the first case, namely (26).

**REMARK 4**. The mathematical transformations carried out highlight a special case (*in my opinion*), in which two sides of the triangle under investigation "intersect" at infinity.

**REMARK 5**. From the obtained identities (8), it becomes evident that the triangles collapse to a point (*in my opinion*) on the plane without "distortions", as a result of which (32) becomes possible:

$$a = b = c = 0 \Rightarrow \left(\left(\frac{0}{2}\right)^2 - \left(\frac{0}{2}\right)^2\right)^2 + (0)^2 = \left(\frac{0}{2} \cdot 0\right)^2 \tag{32}$$

**If** $q = 2$, a contradiction arises. In this case, the number $\left(\frac{a}{2}\right)^2 - \left(\frac{c}{2}\right)^2$ becomes odd, since $a \equiv 0(mod4)$. Consequently, in (21), $q \neq 2$. From this, it also follows that $q \neq 6$. Additionally, according to conditions (11) and (12), $q \neq 4$.

**If** $q = 3$, then, using (4), we write the following:

$$\begin{cases} m_a = \frac{1}{2}\sqrt{2b^2 + 2c^2 - a^2} \\ m_b = \frac{1}{2}\sqrt{2a^2 + 2c^2 - b^2} \\ m_c = \frac{1}{2}\sqrt{2a^2 + 2b^2 - c^2} \end{cases} \Rightarrow \begin{cases} m_a = \frac{1}{2}\sqrt{3b^2 + 3c^2 - (a^2 + b^2 + c^2)} \\ m_b = \frac{1}{2}\sqrt{3a^2 + 3c^2 - (a^2 + b^2 + c^2)} \\ m_c = \frac{1}{2}\sqrt{3a^2 + 3b^2 - (a^2 + b^2 + c^2)} \end{cases} \quad (\mathbf{33})$$

According to condition (33), it is evident:

**If** one median of an *integer triangle* is divisible by 3, then $(a^2 + b^2 + c^2) \equiv 0(mod3)$. This implies that all medians are divisible by 3.

The condition $(a^2 + b^2 + c^2) \equiv 0(mod3)$ is possible only in two cases:

**1**. Either all sides of the triangle are divisible by 3. Then all medians and the area of the *integer triangle* are also divisible by 3. This contradicts condition (9).

**2**. Or none of the sides of the triangle is divisible by 3. This contradicts Corollary 2.

**Consider** a triangle $\Delta ABC$ with integer sides and integer medians.

If all sides of $\Delta ABC$ are divisible by 3, then, in view of (4), we have

$(a, b, c, m_a, m_b, m_c) \geq 3$ for triangle $\Delta ABC$. This contradicts condition (9).

By (4), it is clear that all sides of $\Delta ABC$ are even. Suppose that none of the sides of $\Delta ABC$ is divisible by 3. Then, again by (4), all medians of $\Delta ABC$ are divisible by 3.

By the Lemma, the triangles $\Delta ABC$ and $\Delta OPC$ are linked and exist simultaneously. Consequently, $\frac{2}{3}m_a,\ \frac{2}{3}m_b,\ \frac{2}{3}m_c$ are the sides of triangle $\Delta OPC$, while $\frac{1}{2}a,\ \frac{1}{2}b,\ \frac{1}{2}c$ are its medians. Hence, none of the medians of $\Delta OPC$ is divisible by 3. Therefore, at least one of the sides $\frac{2}{3}m_a,\ \frac{2}{3}m_b,\ \frac{2}{3}m_c$ of triangle $\Delta OPC$ must be divisible by 3.

| | $\Delta ABC$ | $\Delta OPC$ |
|---|---|---|
| Sides: | $316, 262, 254$ | $136, 170, 174$ |
| Medians: | $204, 255, 261$ | $158, 131, 127$ |

**Corollary 2**: Thus, when we consider an arbitrary triangle with integer sides and integer medians under the assumption that none of its sides is divisible by 3, we are simultaneously led to consider another triangle in which exactly one side is necessarily divisible by 3. This, in turn, implies that for one of the equations in (8), we have $q \neq 3$ in (21). Consequently, one of the equations in (8) represents a primitive Pythagorean triple.

**Corollary 3**. There are no Heronian triangles with three integer medians.

**P.S.** Suppose there exist some natural numbers $x$ and $y$ for which the following conditions hold:

$$\begin{cases} x^2 - y^2 = A_{\Delta ABC} \\ 2xy = \left(\frac{a}{2}\right)^2 - \left(\frac{b}{2}\right)^2 \\ x^2 + y^2 = \frac{c}{2} m_c \end{cases}$$

However, we arrive at a contradiction. Previously, we established that only two sides of an *integer triangle* are divisible by $2$ (let these be $c$ and $b$), while the third side (let this be $a$) must be divisible by $4$. This gives

$\left(\frac{a}{2}\right)^2 - \left(\frac{b}{2}\right)^2 \equiv 1(mod2)$. So $2xy \neq \left(\frac{a}{2}\right)^2 - \left(\frac{b}{2}\right)^2$.